
\magnification=\magstep1
\input amstex
\documentstyle{amsppt}
\TagsOnRight
\NoRunningHeads
\newdimen\notespace \notespace=7mm
\newdimen\maxnote  \maxnote=13mm

\define\Rgot{\frak R}
\define\mgot{\frak m}
\define\C {\bold C}
\define\stalk{k[x]_{(x)}}
\define\Spec{\operatorname{Spec}}

\topmatter
\title On the representability of ${\Cal Hilb}^n\stalk$
\endtitle
\author{Roy Mikael Skjelnes}\endauthor
\affil{Department of Mathematics, KTH} \endaffil
\address{R.M. Skjelnes, Dept. of Mathematics, KTH, SE-100 44 Stockholm, Sverige}\endaddress
\email{skjelnes\@math.kth.se}\endemail

\abstract Let $k[x]_{(x)}$ be the polynomial ring $k[x]$ localized in the
 maximal ideal $(x)\subset k[x]$. We study the Hilbert functor parameterizing
 ideals of colength $n$ in this ring {\it  having support at the
 origin}. The main result of this article is that this functor is not
 representable. We also give a complete description of the functor as
 a limit of representable functors
\endabstract
\endtopmatter
\document

\subhead 1. Introduction \endsubhead

Let $k$ be a field. Let $R$ be a local noetherian $k$-algebra with
maximal ideal $P$. The Hilbert functor of $n$-points on $\Spec (R)$,
denoted as ${\Cal Hilb}^n_R$,  is
determined by sending a scheme $T$ to the set
$$\aligned {\Cal Hilb}^n_R(T) \endaligned
= \left\{ \aligned &\text{Closed subschemes }Z\subseteq T\times_k \Spec (R)
  \text{ such that}\\
&\text{the projection } Z \to T \text{ is flat, and where the global}\\
&\text{sections of the fiber }Z_y \text{ is of dimension }n \text{ as a
  }\\
&\kappa (y)\text{-vector space for all points }y \in T.\endaligned \right\}.$$

We let ${\Cal Hilb}^nR(T)\subseteq {\Cal Hilb}^n_R(T)$ be
the set of $T$-valued points $Z$ of ${\Cal Hilb}^n_R$ such that $Z_{red}\subseteq
T\times_k \Spec(R/P)$. Here $Z_{red}$ is the reduced
scheme associated to $Z$. The assignment sending a $k$-scheme $T$ to
the set ${\Cal Hilb}^nR(T)$ determines a contravariant
functor from the category of noetherian $k$-schemes to sets. The
functor ${\Cal Hilb}^nR$ is {\it different} from the Hilbert functor
${\Cal Hilb}^n_R$.

The functor ${\Cal  Hilb}^nR$ with 
$R=\C\{x,y\}$, the ring of convergent power series in two
variables, was introduced by J. Brian\c{c}on in \cite{1}, and its  set of
$\C$-rational points were described. The
motivation behind the present paper was to understand the universal
properties of ${\Cal Hilb}^n\C\{x,y\}$.

Instead of analytic spaces, as considered in \cite{1}, we work
in the category of noetherian $k$-schemes. Primarily our interest were
in the representability of the functor ${\Cal Hilb}^n k[[x,y]]$. However, we
realized that the problems we faced were present for 
${\Cal Hilb}^n\stalk$, where  
$\stalk$ is the local ring of the origin at the line. To illustrate
the difficulties of the representability of ${\Cal
  Hilb}^nk[[x,y]]$ we
will in this paper focus on ${\Cal Hilb}^n\stalk$, the functor parameterizing colength $n$
ideals in $\stalk$, having support in $(x)$.

The scheme $\Spec(k[x]/(x^n))$ is the only closed subscheme of
$\Spec(\stalk)$ whose coordinate ring is of dimension $n$ as a $k$-vector space. It
follows that the functor ${\Cal  Hilb}^n\stalk$ has only one $k$-valued point. Thus in
a naive geometric sense the functor ${\Cal  Hilb}^n\stalk$ is
trivial. We shall see, however, that the functor ${\Cal
  Hilb}^n\stalk$ is not representable! In fact we show in Theorem
(4.8) that ${\Cal Hilb}^n R$ is not representable when $R$ is the
local ring of  a
regular point on a variety.

In addition to Theorem (4.8) which is our main result, we show in Theorem (5.5) that the
  non-representable functor ${\Cal 
  Hilb}^n\stalk$ is pro-represented by $k[[s_1, \ldots, s_n]]$, the
  formal power series ring in $n$-variables. In  Theorem (6.8) we show that
  there exist a natural filtration of ${\Cal  Hilb}^n\stalk$ by representable subfunctors $\{{\Cal
  H}^{n,m}\}_{m\geq 0}$, where  ${\Cal H}^{n,m}$ is a closed subfunctor of ${\Cal
  H}^{n,m+1}$. 

The three theorems (4.8), (5.5) and (6.8) completely describe
${\Cal Hilb}^n\stalk$. The three mentioned 
results are more or less explicit  applications of Theorem (3.5),
which describes the set of elements in ${\Cal Hilb}^n\stalk
(\Spec(A))$ for arbitrary $k$-algebras $A$.

The paper is organized as follows: In Section (2) we recall some results from \cite{2}. In Section (3)
we establish Theorem (3.5). The sections (4), (5) and (6) are applications of Theorem (3.5). In
Section (4) we show that ${\Cal  Hilb}^n\stalk$ is not representable. We
pro-represent ${\Cal  Hilb}^n\stalk$ in Section (5). We give a
filtration of ${\Cal  Hilb}^n\stalk$ by representable subfunctors in
Section (6).

I want to thank  to my thesis advisor Dan Laksov for his help and
assistance during the preparation of the present paper. I thank
Torsten Ekedahl, Trond Gustavsen, Yves Pitteloud and the referee for their comments
and remarks.

\subhead {2. Preliminaries}\endsubhead 

\definition {2.1. Notation} Let $k$ be a field. Let $k[x]$ be the ring
of polynomials in one variable over $k$. The polynomials $f(x)$ in $k[x]$ such that
$f(0)\neq 0$ form a multiplicatively closed subset $S$ in $k[x]$. We write the fraction ring
$k[x ]_S=\stalk$. For every $k$-algebra $A$ we write $A\otimes_k
k[x]=A[x]$. The localization of the $k[x]$-algebra $A[x]$ in the
multiplicatively closed set $S\subset k[x]$ is $A\otimes_k
\stalk$. If $I$ is an ideal in a 
ring $A$ we let $\Rgot(I)$ denote its radical, and if $P$ is a prime
ideal we let $\kappa (P)=A_P/PA_P$ be its residue field.
\enddefinition

\proclaim{Lemma 2.2} Let $A$ be a $k$-algebra. Let $I\subseteq A\otimes_k \stalk$ be an ideal
  such that $A\otimes_k \stalk /I$ is a free $A$-module of rank
  $n$. Then the following two assertions hold:
\roster
\item The classes of $1, x, \ldots ,x^{n-1}$ form an $A$-basis for
  $A\otimes_k \stalk/I$.
\item  The ideal $I$ is generated by a
  unique $F(x)=x^n-u_1x^{n-1}+\dots +(-1)^nu_n$ in $A[x]$.
\endroster
\endproclaim

\demo{Proof} See \cite{3}, Lemma (3.2) for a proof of the first
assertion. The second assertion follows from 
\cite{3}, Theorem (3.3). 
\enddemo

\proclaim{Proposition 2.3} Let $A$ be a $k$-algebra. Let 
  $I\subseteq A\otimes_k \stalk$ be an ideal with residue ring
  $M=A\otimes_k \stalk /I$. Assume that 
\roster
\item There is  an inclusion of ideals $(x)\subseteq \Rgot(I)$ in
  $A\otimes _k \stalk$.
\item The $A$-module $M=A\otimes_k \stalk /I$ is flat.
\item For every prime ideal $P$ in $A$ we have that $M\otimes_A
  \kappa (P)$ is of dimension $n$ as a $\kappa (P)$-vector space.
\endroster
Then $M$ is a free $A$-module of rank $n$.
\endproclaim

\demo{Proof} We first show that $M\otimes_A A_P$ is free for every
  prime ideal $P$ in $A$. Thus we  assume that $A$
  is a local $k$-algebra. Assumption (1) is equivalent to the
  existence of an integer $N$ such that we have an inclusion of ideals 
  $(x^N)\subseteq I$ in $A\otimes_k \stalk$. Consequently we have a
surjection
$$ A\otimes_k \stalk /(x^{N}) \rightarrow M=A\otimes_k \stalk /I.\tag{2.3.1}$$
We have that $A\otimes_k \stalk/(x^{N})=A[x]/(x^{N})$. It follows from the surjection
(2.3.1)  that $M$ is generated by the classes of $1, x, \ldots
,x^{N-1}$. In particular $M$ is finitely generated. A flat and
finitely generated module over a local ring is free, see \cite{4} Theorem
(7.10). Hence by Assumption (2) we have that $M$ is a free
$A$-module. By Assumption (3) we have that the rank of $M$ is
$n$. 

Thus  we have proven that $M\otimes_AA_P$ is free of rank $n$ for every
prime ideal $P$ in $A$. It then follows by  Assertion (1) of
Lemma (2.2) that $M\otimes_A A_P$ has a basis given by the
classes of $1,x, \ldots ,x^{n-1}$. Since the classes of $1, x, \ldots
,x^{n-1}$ form a basis for $M\otimes_AA_P$ for every prime ideal $P$
of $A$, it follows that $1, x,\ldots ,x^{n-1}$ form a basis for $M$. 
\enddemo

\proclaim{Theorem 2.4} Let $A$ be a $k$-algebra and let $F(x)$ in $A[x]$ be a polynomial where $F(x)=x^n-u_1x^{n-1}+\dots +(-1)^nu_n$. The
  following three assertions are equivalent.
\roster
\item  For all maximal ideals $P$ of $A$ with residue map $\varphi
  : A\rightarrow A/P$, the roots of
  $F^{\varphi}(x)=x^n-\varphi(u_1)x^{n-1}+\dots +(-1)^n\varphi (u_n)$
  in the algebraic closure of $A/P$ are zero or transcendental over
  $k$.
\item The ring $A\otimes_k \stalk/(F(x))$ is canonically isomorphic to
  $A[x]/(F(x))$.
\item The $A$-module  $A\otimes_k \stalk/(F(x))$ is free of rank
  $n$ with a basis consisting of the classes of $1,x ,\ldots ,x^{n-1}$.
\endroster
\endproclaim

\demo{Proof} See \cite{3}, Assertions (1), (4)  and (5) of Theorem
(2.3).
\enddemo

\proclaim {Corollary 2.5} Let $F(x)=x^n-u_1x^{n-1}+\dots +(-1)^nu_n $
  be an element of
  $A[x]$. Assume that the
  coefficients $u_1, \ldots ,u_n$  are in the Jacobson
  radical of $A$. Then we have that $M=A\otimes_k
  \stalk/(F(x))$ is
  canonically isomorphic to $A[x]/(F(x))$. In particular we have a
  canonical isomorphism $M= A[x]/(F(x))$ (17) when $A$ is local and the coefficients $u_1,
  \ldots ,u_n$ of $F(x)$ are in the maximal ideal of $A$.
\endproclaim

\demo {Proof} Let $P$ be a maximal ideal, and let $\varphi : A
\rightarrow A/P$ be the residue map.  We have that
$F^{\varphi}(x)=x^n$ since the coefficients $u_1, \ldots ,u_n$
of $F(x)$ are in the Jacobson radical of $A$. Consequently the roots of $F^{\varphi}(x)$ are
zero, and the Assertion (1) of the Theorem is satisfied.
\enddemo

\proclaim{Corollary 2.6} Assume that $F(x)$ in $A[x]$ is such that
  the assertions of the Theorem are satisfied. Then an inclusion of ideals $(x^N)\subseteq (F(x))$ in
  $A\otimes_k \stalk $ is equivalent to an inclusion of ideals
  $(x^N)\subseteq (F(x))$ in $A[x]$.
\endproclaim

\demo{Proof} Obviously an inclusion of ideals in $A[x]$ extends to an
inclusion of ideals in the fraction ring $A\otimes_k
\stalk$. Consequently it suffices to show that an inclusion
$(x^N)\subseteq (F(x))$ in $A\otimes_k \stalk$ gives an inclusion
 $(x^N)\subseteq (F(x))$ in $A[x]$. Assume that we have an
inclusion of ideals $(x^N)\subseteq (F(x))$ in $A\otimes_k \stalk$,
or equivalently  a surjection
$$ A\otimes_k \stalk /(x^N) \rightarrow A\otimes_k \stalk /(F(x)).
\tag{2.6.1}$$
We have that $F(x)$ in $A[x]$ satisfies the conditions in the
Theorem. Hence we have a canonical isomorphism $A\otimes_k
\stalk/(F(x))=A[x]/(F(x))$. Then the   surjection (2.6.1)  gives a surjection
$A[x]/(x^N)\rightarrow A[x]/(F(x))$ which is equivalent to  an
inclusion of ideals $(x^N)\subseteq (F(x))$ in $A[x]$. 
\enddemo

\subhead {3. Polynomials with nilpotent coefficients}\endsubhead

The purpose of this section is to establish Theorem (3.5). Applications of Theorem
(3.5) is given in Sections (4), (5) and (6).
 
\definition{3.1. Set up and Notation} We will study ideals generated by monic
polynomials with nilpotent coefficients. For this purpose we introduce
the following terminology; Let $A$ be a commutative ring, and let $A[t_1,
\ldots ,t_n]$ be the polynomial ring over $A$  in the variables $t_1, \ldots
,t_n$. Let $s_i(t)=s_i(t_1, \ldots ,t_n)$ be the $i$'th elementary
symmetric function in the variables $t_1, \ldots t_n$. The elementary symmetric functions $s_i(t)$ are homogeneous in the
variables $t_1, \ldots ,t_n$, having degree
$\text{deg}(s_i(t))=i$. We let $A_0=A$ and consider the ring of symmetric functions
$A[s_1(t), \ldots ,s_n(t)]=\oplus_{i\geq 0}A_i$ as graded
in $t_1, \ldots ,t_n$. For every positive integer $d$ we have the
ideal $\oplus_{i\geq d}A_i\subseteq A[s_1(t), \ldots ,s_n(t)]$. We
denote the residue ring by
$$Q_d :=A[s_1(t), \ldots ,s_n(t)] /\oplus _{i\geq d}A_i.\tag{3.1.1}$$
\enddefinition

\proclaim{Lemma 3.2} Let  $u_1, \ldots ,u_n$ be nilpotent elements in
  a ring $A$. Then the homomorphism $u :A[s_1(t), \ldots ,s_n(t)]\rightarrow
  A$, determined by $u(s_i)=u_i$ for $i=1, \ldots ,n$, factors
  through $Q_d$ for some integer $d$.
\endproclaim

\demo{Proof} The coefficients $u_1, \ldots ,u_n$ are nilpotent by
assumption.  Hence there exist integers $n_i$ such that $u_i^{n_i}=0$
for every $i=1, \ldots ,n$. Let $\tau =\text{max}\{n_i\}$, and let
$d=\tau+2\tau +\dots +n\tau$. We claim that $u : A[s_1(t), \ldots ,s_n(t)] \rightarrow A$ maps
 $\oplus_{i\geq d}A_i$ to zero. It is enough to show that monomials
 $m(s_1(t), \ldots ,s_n(t))$ of degree $\geq d$ are mapped to zero. We have that $m(s_1(t), \ldots
 ,s_n(t))=s_1(t)^{e_1}s_2(t)^{e_2}\dots s_n(t)^{e_n}$ where $e_1+2e_2
 +\dots + ne_n=\text{deg}(m(s_1(t), \ldots ,s_n(t)))$. It follows that
 at least one $e_j\geq \tau$, and consequently $u_j^{e_j}=0$. Thus we have that $
u (s_1(t)^{e_{1}}s_2(t)^{e_{2}}\dots
s_n(t)^{e_{n}})=u_1^{e_{1}}u_2^{e_{2}}\dots u_n^{e_{n}} =0$. 
\enddemo

\definition{3.3. Polynomials with nilpotent coefficients}  For every monic polynomial $F(x)=x^n-u_1x^{n-1} +\dots +(-1)^nu_n$
in $A[x]$ we let $u_F:A[s_1(t), \ldots ,s_n(t)]\rightarrow A$ be the $A$-algebra
homomorphism determined by $u_F(s_i(t))=u_i$ for $i=1, \ldots ,n$. Let
$$\Delta (t,x)=\prod_{i=1}^n (x-t_i)=x^n-s_1(t)x^{n-1}+\dots +(-1)^ns_n(t).
\tag{3.3.1}$$
If $D(t,x)=D(t_1, \ldots ,t_n,x)$ is symmetric in the variables $t_1,
\ldots ,t_n$, we let $D^{u_F}(x)$ in $A[x]$ be the image of $D(t,x)$ by the  map
$u_F\otimes 1: A[s_1(t), \ldots,s_n(t)][x]\rightarrow A[x]$. In particular we have
that $\Delta ^{u_F}(x)=F(x)$.

For every non-negative integer $p$ we define $d_p(t_i,x)$ in $A[t_1,
\ldots ,t_n,x]$ by
$$d_p(t_i,x)=(x+t_i)(x^2+t_i^2)\dots (x^{2^{p}}+t_i^{2^p}).
\tag{3.3.2}$$
It follows by induction on $p$ that
$(x-t_i)d_p(t_i,x)=x^{2^{p+1}}-t_i^{2^{p+1}}$. We let 
$$ D_p(t,x)=\prod_{i=1}^n d_p(t_i,x).
\tag{3.3.3}$$
For every non-negative integer $N$, we let $s_i(t^N)=s_i(t_1^N,
, \ldots , t_n^N)$, which is a homogeneous symmetric function in the
variables $t_1, \ldots ,t_n$. We have that the degree of $s_i(t^N)$ is
$\text{deg}(s_i(t^N))=iN$. Both $\Delta (t,x)$ and $D_p(t,x)$ are symmetric in the variables
$t_1, \ldots ,t_n$. Their product is
$$\aligned
\Delta(t,x)D_p(t,x) &= \prod _{i=1}^n (x^{2^{p+1}}-t_i^{2^{p+1}}) \\
  &= x^{2^{p+1}n}- s_1(t^{2^{p+1}})x^{2^{p+1}(n-1)} +\dots
  +(-1)s_n(t^{2^{p+1}}).
\endaligned \tag{3.3.4}$$
\enddefinition

\proclaim{Proposition 3.4} Let $F(x)=x^n-u_1x^{n-1}+\dots +(-1)^n
  u_n$ be an element of  $A[x]$. Then the coefficients $u_1,
  \ldots ,u_n$ are nilpotent if and only if we have an inclusion of ideals
  $(x)\subseteq \Rgot(F(x))$ in $A[x]$. 
\endproclaim

\demo {Proof} Assume that the coefficients $u_1, \ldots ,u_n$ of $F(x)$ are nilpotent.We must show that $x\in \Rgot(F(x))$, or equivalently
  that $x^N\in (F(x))$ for some integer $N$.  By Lemma (3.2) the map $u_F:A[s_1(t), \ldots
  ,s_n(t)]\rightarrow A$ determined by $u_F(s_i(t))=u_i$, factors through $Q_d$ for some integer
  $d$. Let $p$
  be an integer such that $2^{p+1}\geq d$. The function $D_p(t,x)$
  (3.3.3) is symmetric in the
  variables $t_1, \ldots ,t_n$. We will show that $D_p^{u_F}(x)$ in $A[x]$ is such
  that $F(x)D_p^{u_F}(x)=x^N$.
  The  product $\Delta (t,x)D_p(t,x)$ is given in (3.3.4), and the
  degree of the symmetric functions $s_i(t^{2^{p+1}})=i2^{p+1}\geq d$. Consequently the class of $\Delta (t,x)D_p(t,x)$
  in $Q_d[x]$ equals $x^{2^{p+1}n}$. We obtain that 
$$
x^{2^{p+1}n} = \Delta^{u_F} (x)D^{u_F}_p(x)=F(x)D_p^{u_F}(x),
\tag{3.4.1}$$
in $A[x]$. Hence we have that $(x)\subseteq \Rgot(F(x))$.

Conversely, assume that we have an inclusion of ideals $(x)\subseteq \Rgot
(F(x))$ in $A[x]$. Then there exist a $G(x)$ in $A[x]$ such that
$x^N=F(x)G(x)$ for some integer $N$. Let $P$ be a prime ideal of $A$,
and let $\varphi : A\rightarrow \kappa(P)=K$ the the residue map. Let
$F^{\varphi}(x)$ and $G^{\varphi}(x)$ be the classes of $F(x)$ and
$G(x)$, respectively, in $K[x]$. We have
$$
x^N=F^{\varphi}(x)G^{\varphi}(x)=(x^n -\varphi(u_1)x^{n-1}+\dots +(-1)^n\varphi(u_n))G^{\varphi}(x),
\tag{3.4.2}$$
in $K[x]$. The ring $K[x]$ is a unique factorization domain, hence  $\varphi(u_i)=0$ for $i=1, \ldots
,n$. Therefore the classes of $u_i$ are  zero in $A/P$ for all prime ideals
$P$ of $A$. We have shown that  $u_1, \ldots ,u_n$ are nilpotent. 
\enddemo

\proclaim{Theorem 3.5} Let $A$ be a $k$-algebra, and let  $I\subseteq
  A\otimes_k \stalk$ be an ideal. Write the residue ring as
  $M=A\otimes_k \stalk/I$. The following two assertions are
  equivalent.
\roster
\item  M is a flat $A$-module such that for every prime ideal $P$
  in $A$ we have that $M\otimes_A\kappa (P)$ is of dimension $n$ as a
  $\kappa (P)$-vector space, and  we have an inclusion of ideals
  $(x)\subseteq \Rgot(I)$ in $A\otimes_k \stalk$.
\item  The ideal $I$ is generated by an element $F(x)$ in $A[x]$,
  of the form $F(x)=x^n-u_1x^{n-1}+\dots +(-1)^nu_n$, where the
  coefficients $u_1, \ldots , u_n$ are nilpotent.
\endroster
\endproclaim

\demo{Proof} Assume that Assertion (1) holds.  By Proposition
(2.3) we have  that  $M$ is a free $A$-module of rank $n$. It
 follows from  Lemma (2.2) that the ideal $I$ is generated by  a unique  $F(x)=x^n-u_1x^{n-1}+\dots +(-1)^nu_n$ in
$A[x]$, and that the classes of $1, x, \ldots ,x^{n-1}$ form a basis
for $M$.  Consequently $F(x)$ in $A[x]$ is such that the
assertions of Theorem (2.4) hold.
By assumption there is an inclusion of ideals $(x)\subseteq \Rgot (F(x))$ in $A\otimes_k
\stalk$. Or equivalently that $(x^N)\subseteq (F(x))$ in $A\otimes_k
\stalk$ for some integer $N$. By Corollary (2.6) we get an inclusion of ideals
$(x^N)\subseteq (F(x))$ in $A[x]$. It follows from
Proposition (3.4) that the coefficients $u_1,\ldots ,u_n$
of $F(x)$ are nilpotent.

Conversely, assume that Assertion (2) holds. Since the coefficients
$u_1, \ldots ,u_n$ of $F(x)$ are nilpotent, we get by Corollary (2.5)
that $F(x)$ is such that the assertions of Theorem (2.4) is
satisfied. Thus $M=A\otimes_k \stalk /(F(x))$ is a free $A$-module of
rank $n$. In particular we have that $M$ is a flat $A$-module such
that $M\otimes_A \kappa (P)$ is of rank $n$, for every prime ideal $P$
in $A$. What is left  to prove is the inclusion of ideals $(x)\subseteq
\Rgot(F(x))$ in $A\otimes_k \stalk$. It follows from Proposition
(3.4)  that there is an inclusion of ideals $(x)\subseteq \Rgot(F(x))$ in
$A[x]$. Consequently there exist an integer $N$ such that we have an
inclusion  $(x^N)\subseteq (F(x))$ in $A[x]$. By  Corollary (2.6) we get an inclusion of
ideals $(x^N)\subseteq (F(x))$ in $A\otimes_k \stalk$. We
have proven the Theorem. 
\enddemo
    
\subhead{4. The non-representability of  ${\Cal Hilb}^n \stalk $}\endsubhead

In this section we define for every local noetherian $k$-algebra $R$,
the functor ${\Cal Hilb}^nR$. We will show in Theorem (4.8)
that the functor ${\Cal Hilb}^nR$ is not representable when
$R$ is the local ring of a regular point on a variety.

\definition{4.1. Notation} If $Z$ is a scheme, we let $Z_{red}$ be the
associated reduced scheme. Given   a morphism of
schemes $Z\rightarrow T$. The fiber over a given 
point $y\in T$ we write as $Z_y=Z\times_T
\Spec(\kappa (y))$. Here $\kappa (y)$ is the residue field of
the point $y\in T$. 
\enddefinition

\proclaim{Lemma 4.2} Let  $I$ and $J$ be two ideals in a ring
  $A$. Assume that $I$ is finitely generated.  Then an inclusion $I\subseteq \Rgot(J)$ is
equivalent to the existence of an integer $N$ such that
$I^N\subseteq J$.
\endproclaim

\demo {Proof} Let  $x_1, \ldots ,x_m$ be a set of generators for the
ideal $I$.  Assume that we have an  inclusion of ideals
$I\subseteq \Rgot(J)$. It follows that there exist
integers $n_i$ such that $x_i^{n_i} \in J$, for $i=1, \ldots ,m$. Thus
we have that $I^N\subseteq J$, when $N\geq \sum_{i=1}^m(n_i-1)+1$.
The converse is immediate, and we have proven the Lemma. 
\enddemo

\proclaim{Lemma 4.3} Let  $I$ be an ideal in a noetherian $k$-algebra
   $R$. Let $T$ be a noetherian $k$-scheme. Suppose that
  $Z\subseteq T\times_k \Spec(R)$ is a closed
  subscheme. Then $Z_{red}\subseteq T\times_k \Spec(R/I)$ if and only if there exist an integer $N=N(Z)$
  such that $Z\subseteq T\times_k \Spec(R/I^N)$.
\endproclaim

\demo{Proof} The scheme $T$ is noetherian and we can find a finite
affine open cover $\{U_i\}$ of $T$. Thus   $\{U_i\times_k
\Spec(R)\}$ is a finite affine open cover of
$T\times_k\Spec(R)$. It follows from the finite covering of
$T\times_k \Spec(R)$ that  it is enough to prove the
statement for each $U_i\times_k \Spec(R)$. Hence we may assume that $T$ is affine.

Let $T=\Spec(A)$, and let the closed subscheme $Z$ be given by the ideal $J\subseteq
A\otimes_k R$. The image of the natural map $A\otimes_k I
\rightarrow A\otimes_k R$, we write as $I_A$. The ring $R$ is
noetherian, hence  $I\subseteq R$ is
finitely generated. Consequently the ideal $I_A\subseteq A\otimes_k R$ is finitely
generated. It follows  from Lemma (4.2) that  $I_A\subseteq \Rgot(J)$ if and only
if $I_A^N \subseteq J$ for some $N$. We have proven the Lemma.
\enddemo

\definition{4.4. Definition}  Let $R$ be a local noetherian $k$-algebra. Let
$P$ be the maximal ideal of $R$. Let
$n$ be a fixed positive integer. We define for any $k$-scheme $T$ the set
$$ \aligned {\Cal Hilb}^nR(T) \endaligned =\left\{ \aligned  &\text{Closed subschemes }
    Z\subseteq T\times_k \Spec (R),\text{ where the} \\
&\text{projection }Z\to T\text{ is flat, such that the global sections}\\
&\text{of the fiber }Z_y \text{ is of dimension } n \text{ as a }\kappa
(y)\text{-vector space} \\
&\text{for all points } y\in T \text{ and such that }Z_{\operatorname{red}}\subseteq T\times_k \Spec
(R/P). \endaligned \right\}. $$
\enddefinition

\proclaim{4.5. Lemma} The assignment  sending a $k$-scheme $T$ to
the set ${\Cal Hilb}^nR(T)$, determines a contravariant functor from
the category of noetherian $k$-schemes to sets.
\endproclaim

\demo{Proof} Let 
$U \rightarrow T$ be a morphism of noetherian $k$-schemes. If $Z$ is
a $T$-valued point of ${\Cal Hilb}^nR$ we must show that $Z_U=U\times_TZ$
is a element of ${\Cal Hilb}^nR(U)$. The only non-trivial part of
the claim
is to show that $Z_U$ is supported at $U\times_k
\Spec(R/P)$. 

Since $Z$ is supported at $T\times_k \Spec(R/P)$ there exist by
Lemma (4.3) an integer $N$ such that $Z\subseteq T\times_k
\Spec(R/P^N)$. It follows that $Z_U \subseteq
U\times_k \Spec(R/P^N)$. Hence by Lemma (4.3) we have that
$Z_U$ is supported at $U\times_k \Spec(R/P)$. We have
proven the claim. 
\enddemo

\remark{Remark} Note that we restrict ourselves to noetherian
$k$-schemes. It  is not clear whether ${\Cal  Hilb}^nR$
is a presheaf of sets on the category of $k$-schemes. 
\endremark

\remark{Remark} When $R=\C\{x,y\}$, the ring of convergent power
series in two variables, the Definition (4.4) gives the functor of
J. Brian\c{c}on \cite{1}. 
\endremark

\proclaim{Lemma 4.6} Let $R$ be a local noetherian $k$-algebra. Let $P$
be the maximal ideal of $R$, and let $\hat R$ be the $P$-adic
completion of $R$. We have that ${\Cal Hilb}^nR$ is canonically
isomorphic to  ${\Cal Hilb}^n{\hat R}$.
\endproclaim

\demo{Proof}   We have that $\hat R$ is a local ring with maximal ideal
$\hat P=P\otimes_ R \hat R$. Furthermore we have for any positive
integer $N$ that $R/P^N =\hat R/\hat P^N$. It follows that for any
$k$-scheme $T$ we have that
$$ T\times_k \Spec(R/P^N) =T\times_k \Spec(\hat R/\hat
P^N). \tag{4.6.1}$$
Thus if $Z$ is an element of ${\Cal Hilb}^nR(T)$ there is
by Lemma (4.3)
an integer $N$ such that $Z$ is a closed subscheme of $T\times _k
\operatorname{Spec}(R/P^N)$. By (4.6.1) it follows that $Z$ is a closed
subscheme of $T\times _k \operatorname{Spec}(\hat R)$ having support in
$T\times _k \operatorname{Spec}(\hat R/\hat P)$. We get that $Z$ is an element of ${\Cal Hilb}^n\hat R(T)$. 
It is clear that a similar argument shows that the converse also holds; any element $Z\in {\Cal
  Hilb}^n\hat R (T)$ is naturally identified as an element of
${\Cal Hilb}^n R(T)$. We have proven the Lemma. 
\enddemo

\proclaim{Lemma 4.7} Let $A$ be a $k$-algebra. Given a nilpotent element
  $\epsilon $ in $A$, such that the smallest integer $j$ where
  $e^{j}=0$ is $j=2^{m+1}$. Then  the smallest
  integer $N$ such that  we have an inclusion of ideals $(x^N)\subseteq
  (x^n-\epsilon x^{n-1})$ in $A\otimes_k \stalk$ is 
  $N=2^{(m+1)}+n-1$.
\endproclaim

\demo{Proof} We first show that we have an inclusion $(x^N)\subseteq
(x^n-\epsilon x^{n-1})$ in $A\otimes_k \stalk$, with $N=2^{m+1}+n-1$. For every  non-negative integer $p$ we let 
$$ d_p(\epsilon, x)=(x+\epsilon)(x^2+\epsilon^2)\dots (x^{2^p}+\epsilon
^{2^{p}}) \quad \text{in}\quad A[x]. \tag{4.7.1}$$
We have that $(x-\epsilon)d_p(\epsilon, x)=x^{2^{(p+1)}}-\epsilon
^{2^{(p+1)}}$ in $A[x]$. Thus when $p\geq m$, we have that
$(x-\epsilon)d_p(\epsilon, x)=x^{2^{p+1}}$ in $A[x]$. It follows
that there
is an inclusion  of ideals
$$(x^{2^{m+1}+n-1})\subseteq (x^n-\epsilon x^{n-1}) \quad \text{in}\quad
A\otimes_k \stalk.
\tag{4.7.2}$$
We need to show that
$2^{m+1}+n-1$ is the smallest integer such that the  inclusion (4.7.2)
in $A\otimes_k \stalk$
holds. 

Let $N+r=2^{m+1}+n-1$, where $r$ is a non-negative integer. Assume
that we have an inclusion of ideals $(x^N)\subseteq (x^n-\epsilon x^{n-1})$ in
$A\otimes_k \stalk $. The element $\epsilon \in A$ is nilpotent, hence
by Corollary (2.5) we have that 
$F(x)=x^n-\epsilon x^{n-1}$ is such that  the assertions of Theorem
(2.4) are satisfied. It follows by Corollary (2.6) that an inclusion of
ideals $(x^N)\subseteq (F(x))$ in $A\otimes_k \stalk$ is equivalent
with an inclusion of ideals $(x^N)\subseteq (F(x))$ in $A[x]$.
Consequently there exist a $G(x)$ in $A[x]$ such that $x^{N}=(x^n-\epsilon
x^{n-1})G(x)$.  Let $d_m(\epsilon,x)$ in $A[x]$ be the polynomial as
defined in (4.7.1). We have that $(x-\epsilon
)d_m(\epsilon, x)=x^{2^{m+1}}$. Hence we get the following identity in $A[x]$;
$$ (x^n-\epsilon x^{n-1})d_m (\epsilon,x)
=x^{2^{m+1}+n-1}=x^Nx^{r}=(x^{n}-\epsilon x^{n-1})G(x)x^{r}.
\tag{4.7.3}$$
The element $x^n$ is not a zero divisor in the ring $A[x]$. It follows that the element $(x^n-\epsilon
x^{n-1})$ is not a zero divisor in $A[x]$. From the identity in
(4.7.3) we obtain the identity $
 (x^n-\epsilon x^{n-1})(d_m(\epsilon,x) -G(x)x^r)=0 \quad \text{in}\quad
A[x]$, which  implies that $d_m(\epsilon,x)=G(x)x^r$ in $A[x]$. The
polynomial $d_m(\epsilon ,x)$ (4.7.1) has a constant term
$\epsilon^{2^{(m+1)}-1}\neq 0$. Consequently $x$ does not divide
$d_m(\epsilon, x)$. Therefore $r=0$, and  $N=2^{m+1}+n-1$ is the smallest integer such that we have an
inclusion of ideals $(x^N)\subseteq (x^n-\epsilon x^{n-1})$ in
$A\otimes_k \stalk$.
\enddemo

\remark{Remark} When $\epsilon (m)=\epsilon$ is as in Lemma (4.7),
we have that the
closed subscheme  $Z_m=\Spec(A\otimes_k \stalk/(x^n-\epsilon x^{n-1}))\subseteq \Spec(A\otimes_k \stalk)$  is a subscheme of
  $\Spec(A)\times_k \Spec(k[x]/(x^N))$ if and only if $N\geq 2^{m+1}+n-1$.
\endremark

\proclaim{Theorem 4.8} Let $R$ be a local noetherian $k$-algebra
  with maximal ideal $P$. Assume that the $P$-adic completion of $R$
  is  $\hat
  R=k[[x_1, \ldots ,x_r]]$, the formal power series ring in $r>0$  variables. Then  we have that the functor ${\Cal  Hilb}^nR$ is not
  representable in the category of noetherian $k$-schemes.
\endproclaim

\demo{Proof} Write $x=x_1, \ldots ,x_r$, and set $\stalk =k[x_1,
\ldots , x_r]_{(x_1, \ldots ,x_r)}$ the localization of the polynomial
  ring $k[x_1, \ldots ,x_r]$ in the maximal ideal $(x_1, \ldots,
  x_r)$. By Lemma (4.6) it suffices to show that ${\Cal Hilb}^n
  \stalk$ is not representable. 

Assume  that ${\Cal Hilb}^n\stalk $ is representable. Let $H$ be the
noetherian $k$-scheme representing the functor ${\Cal  Hilb}^n\stalk$. Let
$U\in {\Cal  Hilb}^n\stalk (H)$ be the universal family. Then in particular
we have that $U_{red}\subseteq H\times_k \Spec(k)$. Hence, by Lemma (4.3)
there exist an integer $N$ such that we have an closed immersion
$$ U \subseteq H\times_k \Spec(k[x_1, \ldots ,x_r]/(x_1, \ldots ,x_r)^N).
\tag{4.8.1}$$
We let $m$ be an integer such that
$2^{m+1}+n-1>N$. Write $A_m=k[u]/(u^{2^{(m+1)}})$. Let $Z_m=\Spec(A_m\otimes_k \stalk /(x_1^n-\epsilon
x_1^{n-1},x_2, \ldots ,x_r))\subseteq \Spec(A_m\otimes_k \stalk)$, where $\epsilon
\in A_m$ is the class of $u$ in $A_m$.  We have that
$Z_m=\Spec(A_m\otimes_k k[x_1]_{(x_1)}/(x_1^n-\epsilon
x_1^{n-1}))$. It follows from
Theorem (3.5) that $Z_m$ is an $A_m$-valued point of the functor ${\Cal 
  Hilb}^n\stalk$. 

By the universality of the pair $(H,U)$ there exist a  morphism
$\Spec(A_m)\rightarrow H$ such that
$Z_m=\Spec(A_m)\times_H U$. It then follows from the closed immersion
in (4.8.1) that $Z_m\subseteq
\Spec(A_m)\times_k \Spec(k[x]/(x_1, \ldots
,x_r)^N)$. However, since $2^{m+1}+n-1>N$ we have by the remark
following Lemma (4.7), that $\operatorname{Spec}(A_m\otimes_k k[x_1]_{(x_1)}/(x_1^n-\epsilon
x_1^{n-1}))$ is {\it not} a subscheme of
$\operatorname{Spec}(A_m)\times_k \Spec(k[x_1]/(x_1^N))$. Hence we get
that $Z_{m}$ can not be a closed subscheme of
$\Spec(A_m)\times_k \Spec(k[x_1]/(x_1,\ldots ,x_r)^N)$. We have thus reached a contradiction and
proven the Theorem. 
\enddemo

\subhead{5. Pro-representing ${\Cal  Hilb}^n\stalk$}\endsubhead

\definition{5.1. Set up} Let $s_1, \ldots ,s_n$
be independent variables over the field $k$. The completion of the
polynomial ring $k[s_1, \ldots ,s_n]$ in the maximal ideal $(s_1,
\ldots ,s_n)$ we write as $R_n=k[[s_1, \ldots ,s_n]]$. We  will show that
$R_n$ pro-represents the functor $\operatorname{Hilb}^n\stalk$. We recall the
basic notions from \cite{5}.
\enddefinition

\definition{5.2. Notation} Let ${\bold C}_k$ be the category where the
objects are local artinian $k$-algebras with residue field $k$, and
where the morphisms are (local) $k$-algebra homomorphism. If $A$ is an
object of ${\bold C}_k$ we say that $A$ is an {\it artin} ring. 

We write $H_n$ for the restriction of the functor ${\Cal 
  Hilb}^n\stalk$ to the category ${\bold C}_k$. Notice that an artin
  ring $A$, that is an element of the category ${\bold C}_k$ has only
  one prime ideal. The residue  field of the only prime ideal of $A$ is
  $k$. The ideal $(x^n)$ is the  only ideal $I$ of $\stalk $ such that the
  residue ring $\stalk /I$ has dimension $n$ as a $k$-vector space. It
  follows that the covariant functor $H_n $ from the category ${\bold
  C}_k$ to sets, maps an artin ring $A$ to the set
$$ \aligned H_n(A) \endaligned = \left\{ \aligned &\text {Ideals
    }I\subseteq A\otimes_k \stalk \text{ such that the residue ring}\\
& M=A\otimes _k \stalk /I \text{ is a flat } A\text{-module, where}\\
& M\otimes_A k=k[x]/(x^n), \text{ and such that there is an}\\
&\text{inclusion of ideals } (x)\subseteq \Rgot(I) \text { in }
A\otimes_k \stalk .\endaligned
\right\}. \tag{5.2.1}$$
\enddefinition

\remark{Remark} Let ${\Cal Hilb}^n_{\stalk}$ denote the usual
Hilbert functor and consider its restriction to the category ${\bold C}_k$. Thus an $A$-valued
point of ${\Cal Hilb}^n_{\stalk}$ is an ideal $I\subseteq A\otimes_k \stalk$ such
that the residue ring $M=A\otimes_k \stalk /I$ is flat over $A$, and
such that $M\otimes_A k=k[x]/(x^{n})$. We shall show that the
restriction of the Hilbert functor ${\Cal Hilb}^n_{\stalk}$ to the
category ${\bold C}_k$ coincides with the functor $H^n$.

Note that an $A$-valued point $M=A\otimes_k\stalk/I$ of ${\Cal
  Hilb}^n_{\stalk}$ is not {\it a priori} finitely generated as an
  module over $A$. However we have the following general result
  (\cite{2}, Theorem (2.4)).

Let $A$ be a local ring with nilpotent radical. Let $M$ be a flat
  $A$-module, and denote the maximal ideal of $A$ with $P$. If
  $\text{dim}_{\kappa (P)}(M\otimes_A \kappa (P))  =\text{dim}_{\kappa
  (Q)}(M\otimes _A \kappa (Q))=n$, for all minimal prime ideals $Q$ in
  $A$. Then $M$ is a free $A$-module of rank $n$.

It follows that when $A$ is an artin ring, and $M=A\otimes_k \stalk/I$
  is an $A$-valued point of ${\Cal
  Hilb}^n_{\stalk}$, then $M$ is free and
of rank $n$ as an $A$-module. It then follows by Lemma (2.2) that the
ideal $I$ is generated by a monic polynomial
  $F(x)=x^n-u_1x^{n-1}+\dots +(-1)^nu_n$ in $A[x]$. Since we have
that $M\otimes_Ak=k[x]/(x^n)$ we get that the coefficients $u_1,
  \ldots ,u_n$ of $F(x)$ are nilpotent. Hence by Theorem (3.5) we have that
$(x)\subseteq \Rgot (I)$ in $A\otimes_k \stalk$. We have shown that
the two functors ${\Cal Hilb}^n\stalk$ and ${\Cal Hilb}^n_{\stalk}$ coincide when
restricted to the category ${\bold C}_k$ of artin rings.
\endremark

\proclaim{Lemma 5.3} Let $A$ be an artin ring. Let $\psi : R_n=k[[s_1, \ldots ,s_n]] \rightarrow A$ be
  a local $k$-algebra homomorphism. Let
  $F^{\psi}_n(x)=x^n-\psi (s_n)x^{n-1}+\dots +(-1)\psi(s_n)$. Then we
  have that $(F^{\psi}_n(x))\subseteq A\otimes_k \stalk$ is an
  $A$-valued point of $H_n$.
\endproclaim

\demo{Proof} Since
  the map  $\psi$ is local we have that $\psi(s_i)$ is  in the maximal
  ideal $\mgot_A$ of $A$, for each  $i=1, \ldots ,n$. The ring $A$ is
  artin. Consequently  $\mgot_A^q=0$ for some integer
  $q$. It follows that the coefficients $\psi(s_1), \ldots
  ,\psi(s_n)$ of $F^{\psi}_n(x)$ are nilpotent. By Theorem  (3.5) we have
  an inclusion of ideals $(x)\subseteq \Rgot(F^{\psi}_n(x))$ in
  $A\otimes_k \stalk$ and  the residue ring $M= A\otimes_k \stalk/(F^{\psi}_n(x))$ is a flat
  $A$-module such that $M\otimes_A k$ is of dimension $n$ as a
  $k$-vector space. Thus we have proven that the ideal
  $(F^{\psi}_n(x))\subseteq A\otimes_k \stalk$ is an   element  of
  $H_n(A)$. 
\enddemo

\definition{5.4. The pro-couple $(R_n,\xi)$} Let $\mgot$ be
the maximal ideal of $R_n=k[[s_1, \ldots ,s_n]]$. For every positive
integer $q$ we let $s_{q,1}, \ldots ,s_{q,n}$ be the classes of $s_1,
\ldots ,s_n$ in $ R_n/\mgot^q$. It follows from
Lemma
(5.3)  that the ideal generated by $F^{q}_n(x)=x^n-s_{q,1}x^{n-1}+\dots
+(-1)^n s_{q,n}$ in $R/\mgot^q [x]$ generates  an $R_n/\mgot^q$-point
of $H_n$. We get a  sequence 
$$ \xi=\{(F^{q}_n(x))\}_{q\geq 0},\tag{5.4.1}$$
where $(F^{q}_n(x))$ is an $R_n/\mgot^q$-point for every
non-negative integer $q$. Clearly $\xi$ defines a point  in  the projective limit
$\varprojlim_{q}\{H^n(R_n/\mgot^q)\}$. Thus we have that $(R_n,\xi)$ is a
{\it pro-couple} of $H_n$. 

We let $h_R$ be
  the covariant functor from ${\bold C}_k$ to sets, which sends an artin
  ring $A$ to the set of local $k$-algebra homomorphisms
  $\text{Hom}_{k\text{-loc}}(R_n,A)$. We note that a local
  $k$-algebra homomorphism $\psi : R_n \rightarrow A$ factors
  through $R_n/\mgot^q$ for high enough $q$. We get that the  pro-couple $(R_n,\xi)$ induces a morphism of functors $F_{\xi} : h_R
\rightarrow H_n$ which for any artin ring $A$, maps an element $\psi
  \in h_R(A)$ to the element $(F^{\psi}_n(x))$  in $H_n(A)$. Here
  $F^{\psi}_n(x)$ is as in Lemma (5.3). 
\enddefinition

\proclaim{Theorem 5.5} Let $R_n=k[[s_1, \ldots , s_n]]$, and
  let $\xi$ be as in {\rm (5.4.1)}. The morphism of functors $F_{\xi}:h_R \rightarrow H_n$ induced by the pro-couple $(R_n,\xi)$, 
  is an isomorphism.
\endproclaim

\demo{Proof} We must construct an inverse to the morphism $F_{\xi}: h_R
\rightarrow H_n$. Let $A$ be an artin ring, and let $I\subseteq
A\otimes_k \stalk$ be an ideal satisfying the properties of
(5.2.1). We  have that  Assertion (1) of  Theorem
(3.5) holds. Consequently the ideal  $I\subseteq
A\otimes_k \stalk$ is generated by a
unique $F(x)=x^n-u_1x^{n-1}+\dots  +(-1)u_n$ in $A[x]$, where $u_1, \ldots ,u_n$
are nilpotent. The coefficients $u_1, \ldots ,u_n$ of $F(x)$ are in
the maximal ideal of $A$, hence the map $\psi: k[[s_1, \ldots , s_n]] \rightarrow
A$ sending $s_i$ to $u_i$, determines a local $k$-algebra homomorphism. We have thus constructed a morphism of functors $G : H_n
\rightarrow h_R$. It is clear that $G$ is the inverse of $F_{\xi}$. 
\enddemo

\subhead{6. A filtration  of  ${\Cal  Hilb}^n \stalk$ by schemes}\endsubhead

We will in Section (6) show that there is a natural filtration of
${\Cal Hilb}^n\stalk$ by representable functors $\{{\Cal
  H}^{n,m}\}_{m\geq 0}$, where ${\Cal H}^{n,m}$ is a closed subfunctor
of ${\Cal H}^{n,m+1}$ for all $m$. The functors ${\Cal H}^{n,m}$ are
the Hilbert functors parameterizing closed subschemes of length $n$ of
$\Spec (k[x]/(x^{n+m}))$.

An outline of Section (6) is as follows. We will define  the functors ${\Cal H}^{n,m}$ from the
category of $k$-schemes, not necessarily noetherian schemes, to
sets. We then construct schemes $\Spec (H_{n,m})$ which we show
represent ${\Cal H}^{n,m}$. Thereafter we restrict ${\Cal H}^{n,m}$ to
the category on noetherian $k$-schemes, and show that we get an
filtration of ${\Cal Hilb}^n\stalk$.

\definition{6.1. Definition} Let $n>0,m\geq 0$ be integers. In the
  polynomial ring $k[x]$ we have the ideal $(x^{n+m})$ and we  denote the residue ring as
  $R=k[x]/(x^{n+m})=\stalk /(x^{n+m})$. We denote
  by ${\Cal H}^{n,m}={\Cal Hilb}^n_R$ the local Hilbert functor of
  $n$-points on $\Spec{R}$. Thus ${\Cal H}^{n,m}$ is  the
  contravariant functor from the category of
  $k$-schemes to sets, determined by sending a $k$-scheme $T$ to the
  set 
$$ \aligned {\Cal H}^{n,m}(T) \endaligned =\left\{ \aligned  &\text{Closed subschemes }
    Z\subseteq T\times_k \Spec (k[x]/(x^{n+m}) ), \\
&\text{such that the projection }Z\to T\text{ is flat, and where }
\\
&\text{the global sections of the fiber }Z_y\text{ is of dimension } n \\
&\text{as a }\kappa (y)\text{-vector space, for all points }y \in T. \endaligned \right\}. $$\enddefinition

\definition{6.2. Construction of the rings $H_{n,m}$} Let $P_n=k[s_1,
    \ldots ,s_n]$ be the polynomial ring in the variables
    $s_1\ldots ,s_n$ over $k$. Let $m$ be a fixed non-negative integer, and let $y_1, \ldots ,y_m,x$ be algebraic
independent variables over $P_n$. We define
    $F_n(x)=x^n-s_1x^{n-1}+\dots +(-1)^n s_n$ in $P_n[x]$, and we let
$Y_m(x)=x^m+y_1x^{m-1}+\dots +y_m$. The product
$F_n(x)Y_m(x)$ is
$$F_n(x)Y_m(x)=x^{n+m}+C_{m,1}(y)x^{n+m-1}+\dots +
C_{m,n+m}(y).
\tag{6.2.1}$$
As a convention we let
$s_0=y_0=1$, and $y_j=0$ for negative values of $j$. The coefficient $C_{m,i}(y)$
is the  sum of products $(-1)^{j}s_{j}y_{i-j}$, where $j=0,
\ldots ,n$, and $i-j=0, 1, \ldots , m$. We have 
$$\aligned
& C_{m,i}(y) =y_i-s_1y_{i-1}+\dots + (-1)^ns_ny_{i-n} \\
& C_{m,m+j}(y)=(-1)^jy_ms_j+\dots +(-1)^ny_{m+j-n}s_n
\endaligned \quad \aligned
&\text{when } i=1,\ldots ,m.  \\
&\text{when } j=1, \ldots ,n.
\endaligned \tag{6.2.2}
$$
For every non-negative integer $m$ we let
  $I_m\subseteq P_n[y_1, \ldots ,y_m]$ be the ideal generated by the coefficients
 $C_{m,1}(y), \ldots, C_{m,m+n}(y)$. We write
$$ H_{n,m}=P_n[y_1, \ldots ,y_m]/I_m=P_n[y_1, \ldots ,y_m]/(C_{m,1}(y),
\ldots , C_{m,m+n}).
\tag{6.2.3}$$
Using (6.2.2) we note that $C_{m,m}(y) = y_{m}+ C_{m-1,m}(y)$. For every positive integer $m$ we define the $P_n$-algebra
homomorphism 
$$c_m : P_n[y_1, \ldots , y_m] \rightarrow P_n[y_1, \ldots ,y_{m-1}]
\tag{6.2.4}$$
by sending $y_i$ to $y_i$ when $i=1, \ldots ,m-1$, and $y_m$ to
$-C_{m-1,m}(y)$.
\enddefinition

\proclaim{Lemma 6.3} For every non-negative integer $m$ we have
  that the natural map  $P_n \rightarrow P_n[y_1, \ldots
  ,y_m]/(C_{m,1}(y),\ldots, C_{m,m}(y))$ is
  an isomorphism. In particular we get that the map $P_n \rightarrow
  H_{n,m}$ is surjective.
\endproclaim

\demo{Proof}  Consider the homomorphism $c_m$ as defined in (6.2.4). It is clear that $c_m$ is surjective
and that  we get an induced isomorphism 
$$P_n[y_1, \ldots
,y_{m}]/(C_{m,m}(y))\backsimeq P_n[y_1, \ldots ,y_{m-1}].
\tag{6.3.1}$$
When $i\leq m $ we have that $C_{m,i}(y)$ is a function in the
variables $y_1, \ldots ,y_i$. Hence when $i=1, \ldots , m-1$ the
elements $C_{m,i}(y)$ are invariant under the action of $c_{m}$. From (6.2.2) we get that $C_{m,i}(y)=C_{m-1,i}(y)$ when $i=1, \ldots
,m-1$.  It follows by successive use of (6.3.1) that  we get an induced
isomorphism 
$$ P_n[y_1, \ldots ,y_{m}]/(C_{m,1}(y), \ldots
  ,C_{m,m}(y))\backsimeq P_n.
\tag{6.3.2}$$
It is easy to see that the map (6.3.2) composed with the natural map
induced by
$P_n \rightarrow P_n[y_1, \ldots ,y_m]$, is the identity map on
$P_n$. We have proven the Lemma. 
\enddemo

\proclaim{Lemma 6.4} For every positive integer $m$, the
  $P_n$-algebra homomorphism $c_m$ {\rm (6.2.4)} induces
  a surjective  map $H_{n,m} \rightarrow H_{n,m-1}$.
\endproclaim

\demo{Proof} Let $\hat c_m$ be the composite of the residue
map $P_n[y_1, \ldots ,y_{m-1}] \rightarrow H_{n,m-1}$ and  $c_m$. We first show
that we get an induced map $H_{n,m}\rightarrow H_{n,m-1}$. That is, we
show that the ideal $I_m\subseteq P_n[y_1, \ldots ,y_m]$ defining
$H_{n,m}$, is in the kernel of $\hat c_m$.

The ideal $I_m$ is generated by $C_{m,1}(y), \ldots ,C_{m,m+n}(y)$. As
noted in the proof of Lemma (6.3) the elements $C_{m,i}(y)$ are mapped
to $C_{m-1,i}(y)$ when $i=1, \ldots , m-1$, whereas $C_{m,m}(y)$ is in
the kernel of $c_m$. Consequently we need to show that the elements
$C_{m,m+j}(y)$ are mapped to zero by $\hat c_{m}$. Using (6.2.2)
we get that
$$\aligned
C_{m,m+j}(y) &=(-1)^jy_{m}s_j +(-1)^{j+1}y_{m-1}s_{j+1}\dots
 +(-1)^ny_{m+j-n}s_n \\
 &= (-1)^ny_{m}s_j+C_{m-1,m+j}(y) \quad \text{ when } j\leq n-1.
\endaligned \tag{6.4.1}$$
It follows  that $C_{m,m+j}(y)$, for $j=1, \ldots ,n-1$ are
mapped to zero by $\hat {c}_{m}$. The last generator of $I_{m}$ is
$C_{m,m+n}(y)=(-1)^{n}y_{m}s_n$, clearly in the kernel of $\hat
{c}_{m}$. Thus we have proven that the
ideal $I_{m}$ is in the kernel of $\hat c_{m} : P_n[y_1, \ldots ,y_{m}]
\rightarrow H_{n,m-1}$.

We need to show that the induced map $H_{n,m}
\rightarrow H_{n,m-1}$ is surjective. From Lemma (6.3) we have that the
natural map $P_n \rightarrow H_{n,m}$ is surjective for all $m$. Since
the map $c_{m}$ is $P_n$-linear, it follows that the induced map $H_{n,m}
\rightarrow H_{n,m-1}$ is $P_n$-linear and the result follows. 
\enddemo

\definition{6.5. Definition} The natural map $P_n=k[s_1, \ldots ,s_n] \rightarrow
H_{n,m}$ is surjective by Lemma (6.3), for all $m$. We let $s_{m,i}$ be the class of $s_i$ in
$H_{n,m}$, for $i=1, \ldots , n$. Define
$$
F_{n,m}(x)=x^n-s_{m,1}x^{n-1} +\dots +(-1)^n
s_{m,n} \quad \text{in } H_{n,m}[x].
\tag{6.5.1}$$
\enddefinition 

\proclaim{Lemma 6.6} Let $A$ be a $k$-algebra. Given an ideal
  $I\subseteq A\otimes_k \stalk$ such  that the residue ring
  $A\otimes_k \stalk /I$ is  a free
  $A$-module of rank $n$, and such that there is an inclusion of ideals
  $(x^{n+m})\subseteq I$ in $A\otimes_k \stalk $. Then there is a unique
  $k$-algebra homomorphism $\psi: H_{n,m} \rightarrow A$ such that
$$
  F_{n,m}^{\psi}(x)=x^n-\psi(s_{m,1}) x^{n-1}+\dots +(-1)^n
  \psi (s_{m,n})$$
in $A[x]$ generates $I$.
\endproclaim

\demo {Proof} It follows by Assertion (2) of 
Lemma (2.2) that $I$ is generated by a unique
$F(x)=x^n-u_1x^{n-1}+\dots +(-1)^nu_n$ in $A[x]$. By Assertion (1) of Lemma (2.2) the classes of $1,x,  \ldots, x^{n-1}$
form a basis for $M$. Consequently $F(x)$ in $A[x]$ satisfies the
assertions  of
Theorem (2.4). By Corollary  (2.6)  the inclusion of ideals
$(x^{n+m})\subseteq (F(x))$ in $A\otimes_k \stalk$ is equivalent with
the existence of $G(x)$ in $A[x]$ such that $x^{n+m}=F(x)G(x)$. Let
$G(x)=x^m+g_1x^{m-1}+\dots +g_m$ in $A[x]$. The coefficients $g_1,
\ldots ,g_m$ are
uniquely determined by $G(x)$, hence uniquely determined by the ideal
$I$. Let $y_1, \ldots ,y_m$ be independent variables over $k$. We get a well-defined $k$-algebra
homomorphism $\theta :k[s_1, \ldots ,s_n,y_1,
\ldots , y_m]\rightarrow A$ determined by $\theta (s_i)=u_i$ where $i=1,
\ldots ,n$, and
$\theta (y_j)=g_j$ where $j=1, \ldots ,m$.  We have thus constructed a
$k$-algebra homomorphism $\theta : P_n[y_1, \ldots ,y_m] \rightarrow
A$. We will next show that the map $\theta $ factors through
$H_{n,m}$. We have that 
$$\aligned
x^{n+m} &= F(x)G(x)\\
  &=(x^n-u_1x^{n-1}+\dots + (-1)^nu_n)(x^m+g_1x^{m-1}+\dots +g_m)\\ 
  &= x^{n+m}+c_1x^{n+m-1}+ \dots +c_{n+m}
\endaligned \tag{6.6.1}$$
in $A[x]$. It follows that the coefficients $c_j$ where $j=1, \ldots ,m+n$ are
zero in $A$. The homomorphism  $\theta $ induces a map  $P_n[y_1,
\ldots ,y_m][x] \rightarrow A[x]$ which sends $F_{n,m}(x)$ to $F(x)$
and $Y_m(x)=x^m+y_1x^{m-1}+\dots +y_m$ to $G(x)$. It follows that  the
coefficient equations $C_{m,j}(y)$ (6.2.1) where $j=1, \ldots ,m+n$,  are mapped to
$c_j=0$. Hence the homomorphism $\theta :P_n[y_1, \ldots ,y_m]
\rightarrow A$ factors through $H_{n,m}$. Let $\psi : H_{n,m}
\rightarrow A$ be the induced map. We have for each $i=1, \ldots , n$
that  $\psi (s_{m,i})=\theta (s_i)=u_i$. Consequently we get that
$F^{\psi}_{n,m}(x)=F(x)$. We have thus proven the existence of a map
$\psi : H_{n,m} \rightarrow A$ such that $F^{\psi}_{n,m}(x)$ generates
the ideal $I$ in $A\otimes_k \stalk$.

We need to show that the map $\psi $ is the only map with the property
that $(F^{\psi}_{n,m}(x))=I$. Let $\psi' : H_{n,m} \rightarrow A$ be a
$k$-algebra homomorphism such that $F^{\psi '}_{n,m}(x)$ generates the
ideal $I$ in $A\otimes_k \stalk$. By Assertion (2) of Lemma (2.2) the ideal
$I\subseteq A\otimes_k \stalk$ is generated by a unique monic
polynomial $F(x)$ in $A[x]$. It follows that we must have
$F^{\psi '}_{n,m}(x)=F(x)$.  Thus if $u_1, \ldots, u_n$ are the
coefficients of $F(x)$, we get that $\psi '(s_{m,i})=u_i$. A
$k$-algebra homomorphism $H_{n,m}\rightarrow A$
is determined by its action on  $s_{m,1}, \ldots ,s_{m,n}$. Hence
$\psi =\psi'$. We have proven
the Lemma. 
\enddemo

\proclaim{Proposition 6.7} The functor ${\Cal H}^{n,m}$
  is represented by $\Spec(H_{n,m})$ {\rm (6.2.3)}. The
    universal family is given by $\Spec(H_{n,m}[x]/(F_{n,m}(x)))$.
\endproclaim

\demo {Proof} We first show that $\Spec(H_{n,m}[x]/(F_{n,m}(x)))$ is
an $H_{n,m}$-valued point of ${\Cal H}^{n,m}$. We have that $F_{n,m}(x)=x^n-s_{m,1}x^{n-1}+\dots
+(-1)^ns_{m,n}$ in $H_{n,m}[x]$. Since $F_{n,m}(x)$ is of degree $n$
and has leading
coefficient 1, we have that $H_{n,m}[x]/(F_{n,m}(x))$ is a free
$H_{n,m}$-module of rank $n$. By the identity in (6.2.1) and the
construction of $H_{n,m}$  we  have an inclusion of
ideals $(x^{n+m})\subseteq (F_{n,m}(x)) $ in $H_{n,m}[x]$. Thus we
have that $H_{n,m}[x]/(F_{n,m}(x))=H_{n,m}\otimes_k R/(F_{n,m}(x))$,
where $R=k[x]/(x^{n+m})$, and consequently $\Spec(H_{n,m}[x]/(F_{n,m}(x)))$ is an
$H_{n,m}$-valued point of ${\Cal H}^{n,m}$. 

We then  have a morphism of
  functors $F:\text{Hom}(-, \Spec(H_{n,m}))\rightarrow {\Cal
    H}^{n,m}$, which we claim is an isomorphism. 

Let $T$ be a $k$-scheme and let $Z$ be an $T$-valued point of ${\Cal
  H}^{n,m}$.  Let
$p:T\times_k \Spec(k[x]/(x^{n+m})) \rightarrow T$ be the projection on
the first factor. Let $\Spec(A)=U\subseteq T$ be an open affine
  subscheme and let the
closed subscheme $Z\cap p^{-1}(U)\subseteq U\times_k
\Spec(k[x]/(x^{n+m}))$ be given by the ideal $J\subseteq A\otimes_k
k[x]/(x^{n+m})$. Let
$I$ be the inverse image of $J$ under the residue map $A\otimes_k \stalk \to
A\otimes_k k[x]/(x^{n+m})$. 

It follows from the definition of the functor ${\Cal H}^{n,m}$ that the ideal $I$ satisfies the conditions of
Proposition (2.3). Hence $A\otimes_k \stalk /I$ is
a free $A$-module of rank $n$. We have by definition an inclusion of
ideals $(x^{n+m})\subseteq I$ in $A\otimes_k\stalk$. Consequently
we get by Lemma (6.6)  a unique
map $f_U:U\rightarrow \Spec(H_{n,m})$ such that $Z\cap p^{-1}(U)=U\times
_{H_{n,m}}\Spec(H_{n,m}[x]/(F_{n,m}(x)))$. 

Thus, if $\{U_i\}$ is an open
  affine  covering of $T$, we get maps $f_i : U_i \rightarrow
  \Spec(H_{n,m})$ with the property that 
$$ Z\cap p^{-1}(U_i)= U_i\times_{H_{n,m}}\Spec(H_{n,m}[x]/(F_{n,m}(x))).
\tag{6.7.1}$$
The maps $f_i:U_i \rightarrow \Spec(H_{n,m})$ are unique with
  respect to the property (6.7.1). Hence the maps $f_i$ glues together 
   to  a unique map $f_Z:
  T\rightarrow \Spec(H_{n,m})$ such that
  $Z=T\times_{H_{n,m}}\Spec(H_{n,m}[x]/(F_{n,m}(x)))$. It follows from the uniqueness of the map
  $f_Z$ that the assignment sending a $T$-valued point
  $Z$ to the morphism $f_Z$  puts up
  an bijection between the set ${\Cal H}^{n,m}(T)$ and the set
  $\text{Hom}(T,\Spec(H_{n,m}))$. We have proven the
  Proposition.
\enddemo

\proclaim {Theorem 6.8} Let $n$ be a fixed positive integer. There is a filtration of the functor ${\Cal  Hilb}^n\stalk$ by an
  ascending chain  of representable functors
$$ {\Cal H}^{n,0}\subseteq {\Cal H}^{n,1}\subseteq
  {\Cal H}^{n,2}\subseteq \dots ,
$$
where ${\Cal H}^{n,m}$ is a closed subfunctor of ${\Cal
  H}^{n,m+1}$, for every $m$.
\endproclaim

\demo {Proof} By Proposition (6.7) the functors ${\Cal H}^{n,m}$ are
represented by $\Spec(H_{n,m})$ where the  universal family is
given by
$U_{n,m}=\Spec(H_{n,m}[x] /(F_{n,m}(x))$. Let $c_{m+1} : H_{n,m+1}
\rightarrow H_{n,m}$ be the surjective map of Lemma (6.3). It follows
from the $P_n$-linearity of $c_{m+1}$ that the induced map
$H_{n,m+1}[x] \rightarrow H_{n,m}[x]$ maps $F_{n,m+1}(x)$   to
$F_{n,m}(x)$. Consequently we have that  $\Spec(H_{n,m})$ is a
closed subscheme of $\Spec(H_{n,m+1})$ such that $U_{n,m+1}\times
_{H_{n,m+1}}\Spec(H_{n,m})=U_{n,m}$.  Hence we have that ${\Cal H}^{n,m}$ is a closed
subfunctor of ${\Cal H}^{n,m+1}$.

From the constructions (6.2.3)
of the rings $H_{n,m}$ it is evident that they are noetherian. It
follows that the restriction of the functor ${\Cal H}^{n,m}$ to
the category of {\it  noetherian} $k$-schemes, is represented by
$\Spec(H_{n,m})$.

That the functors
$\{{\Cal H}^{n,m}\}_{m\geq 0}$ give a filtration of the functor
${\Cal 
  Hilb}^n\stalk$, follows from Lemma (4.3). Indeed, given an
noetherian $k$-scheme $T$ and  let $Z$ be a 
$T$-valued point of ${\Cal  Hilb}^n\stalk$.  Then there exist an integer
$N=N(Z)$ such that $Z\subseteq T\times_k
\Spec(k[x]/(x^N))$. Consequently the $T$-valued point $Z$ of
${\Cal  Hilb}^n\stalk$ is a $T$-valued point of ${\Cal H}^{n,N-n}$.  We
have proven the Theorem. 
\enddemo

\example{6.9 Examples of $H_{n,m}$} The rings $H_{n,m}$ are all of
the form $k[s_1, \ldots ,s_n]/J_m$, where $J_m$ is generated by $n$
elements.  With  $n=1$ it is not
difficult to solve the equations (6.2.2). We get that $H_{1,m}
=k[u]/(u^{m+1})$. Thus we have that the scheme $\Spec{k[x]/(x^{m+1})}$
itself represents the Hilbert functor ${\Cal H}^{1,m}$ of 1-point on
$\Spec(k[x]/(x^{m+1})$, for all non-negative integers $m$.

In general, with $n>1$ a  description of the generators of
the ideal $J_m$ is not known, even though they can be
recursively solved. For instance, we have
$$\aligned
&H_{2,1}=k[x,y]/(x^2,xy)  \\
&H_{2,2}=k[x,y]/(x^3-2xy,x^2y-y^2) \\
&H_{2,3}=k[x,y]/(x^4-3x^2y +y^2, x^3y-2xy^2).
\endaligned $$\endexample

\enddocument
\end